\documentclass[11pt]{article}
\usepackage{amsmath}
\usepackage{amssymb}

\usepackage{theorem}
\numberwithin{equation}{section}

\newtheorem{theorem}{Theorem}[section]
\newtheorem{proposition}[theorem]{Proposition}

\newtheorem{lemma}[theorem]{Lemma}

\begin{document}
\title{Improved almost Morawetz estimates for the cubic nonlinear Schr{\"o}dinger equation}
\date{\today}
\author{Ben Dodson}
\maketitle

\noindent \textbf{Abstract:} We prove global well-posedness for the cubic, defocusing, nonlinear Schr{\"o}dinger equation on $\mathbf{R}^{2}$ with data $u_{0} \in H^{s}(\mathbf{R}^{2})$, $s > 1/4$. We accomplish this by improving the almost Morawetz estimates in \cite{CR}.

\section{Introduction} The cubic, defocusing, nonlinear Schr{\"o}dinger equation on $\mathbf{R}^{2}$,

\begin{equation}\label{1.1}
\aligned
i u_{t} + \Delta u &= |u|^{2} u, \\
u(0,x) &= u_{0}(x) \in H^{s}(\mathbf{R}^{2}),
\endaligned
\end{equation}

\noindent has been the subject of a great deal of research in recent years. It was proved in \cite{CaWe1} that for any $s > 0$, $(\ref{1.1})$ has a local solution on some interval $[0, T]$, $T(\| u_{0} \|_{H^{s}}) > 0$. Moreover, for a solution to fail to extend to a global solution, but instead exist only on a maximal interval $[0, T_{\ast})$,

\begin{equation}\label{1.1.1}
\lim_{t \rightarrow T_{\ast}} \| u(t) \|_{H^{s}(\mathbf{R}^{3})} = \infty.
\end{equation}

\noindent The first progress to proving the existence of a global solution was proved in \cite{CaWe}.

\begin{theorem}\label{t1.1}
$(\ref{1.1})$ has a global solution for $u_{0} \in H^{1}(\mathbf{R}^{2})$.
\end{theorem}

\noindent \emph{Sketch of Proof:} $(\ref{1.1})$ has the conserved quantities

\begin{equation}\label{1.2}
M(u(t)) = \int |u(t,x)|^{2} dx = M(u(0)).
\end{equation}

\begin{equation}\label{1.3}
E(u(t)) = \frac{1}{2} \int |\nabla u(t,x)|^{2} dx + \frac{1}{4} \int |u(t,x)|^{4} dx = E(u(0)).
\end{equation}

\noindent Combining this fact with the Sobolev embedding theorem implies $E(u(0)) \lesssim \| u_{0} \|_{H^{1}(\mathbf{R}^{2})}^{2}$. Since $(\ref{1.3})$ is positive definite, this implies $\| u(t) \|_{H^{1}}^{2} \leq C(\| u_{0} \|_{2}) \| u_{0} \|_{H^{1}}^{2}$ for all time. Thus there exists a solution for all time. $\Box$\vspace{5mm}

\noindent The reader will notice there is a gap between the regularity necessary to prove local well-posedness $(s > 0)$, \cite{CaWe1} and the regularity needed in Theorem $\ref{t1.1}$ to prove a global solution, \cite{CaWe}. Many have undertaken to close this gap. The first progress was made in \cite{B}.

\begin{theorem}\label{t1.2}
If $u_{0} \in H^{s}(\mathbf{R}^{2})$, $s > 3/5$, then $(\ref{1.1})$ has a global solution of the form

\begin{equation}\label{1.4}
\aligned
u(t,x) &= e^{it \Delta} u_{0} + w(t,x), \\
w(t,x) &\in H^{1}(\mathbf{R}^{2}).
\endaligned
\end{equation}
\end{theorem}

\noindent In this case the method of proof was the Fourier truncation method. Take $\phi(\xi) \in C_{0}^{\infty}$,
$$\phi(\xi) = \left\{
                \begin{array}{ll}
                  1, & \hbox{$|\xi| \leq 1$;} \\
                  0, & \hbox{$|\xi| > 2$.}
                \end{array}
              \right.$$

\noindent Then split the initial data into low frequency and high frequency components. $$\hat{u}_{0}(\xi) = \phi(\frac{\xi}{N}) \hat{u}_{0}(\xi) + (1 - \phi(\frac{\xi}{N}) \hat{u}_{0}(\xi) = \hat{u}_{l}(\xi) + \hat{u}_{h}(\xi).$$ Since $\| u_{l} \|_{H^{1}} \lesssim N^{1 - s} \| u_{0} \|_{H^{s}}$, the equation

\begin{equation}\label{1.5}
\aligned
i v_{t} + \Delta v &= |v|^{2} v, \\
v(0,x) &= u_{l},
\endaligned
\end{equation}

\noindent has a global solution with $$E(v(t,x)) \lesssim N^{2 - 2s} \| u_{0} \|_{H^{s}(\mathbf{R}^{2})}.$$ Also, if $s > 3/5$, the equation

\begin{equation}\label{1.6}
\aligned
i w_{t} + \Delta w &= |v + w|^{2}(v + w) - |v|^{2} v, \\
w(0,x) &= u_{h},
\endaligned
\end{equation}

\noindent has a solution on $[0, T]$ of the form $$e^{it \Delta} u_{h} + q(t,x),$$ $$q(t,x) \in H^{1}(\mathbf{R}^{2}) \hspace{5mm} \forall t.$$

\noindent This approach was modified in \cite{CKSTT1} to produce the I-method. The I-operator,

\begin{equation}\label{1.7}
I_{N} : H^{s}(\mathbf{R}^{2}) \rightarrow H^{1}(\mathbf{R}^{2}),
\end{equation}

\noindent is the smooth, radial Fourier multiplier

\begin{equation}\label{1.8}
\aligned
\widehat{I_{N} f}(\xi) &= m_{N}(\xi) \hat{f}(\xi), \\
m_{N}(\xi) &= \left\{
               \begin{array}{ll}
                 1, & \hbox{for $|\xi| \leq N$;} \\
                 (\frac{|\xi|}{N})^{s - 1}, & \hbox{when $|\xi| > 2N$.}
               \end{array}
             \right.
\endaligned
\end{equation}

\noindent From this point on, we will understand that $I$ refers to the I - operator $I_{N}$.

\begin{equation}\label{1.8.1}
\aligned
\| If \|_{H^{1}(\mathbf{R}^{3})} &\lesssim N^{1 - s} \| f \|_{H^{s}(\mathbf{R}^{3})}, \\
\| f \|_{H^{s}(\mathbf{R}^{2})} &\lesssim \| If \|_{H^{1}(\mathbf{R}^{2})}.
\endaligned
\end{equation}

\noindent Therefore, if $$E(Iu(t)) = \frac{1}{2} \int |\nabla Iu(t,x)|^{2} dx + \frac{1}{4} \int |Iu(t,x)|^{4} dx,$$ was a conserved quantity then the existence of a global solution would follow for any $s > 0$. This is not the case, however. Instead, it was proved in \cite{CKSTT1} that

\begin{lemma}\label{l1.2.1}
If $E(Iu(t)) \leq 1$, then there exists $\delta > 0$ such that

\begin{equation}\label{1.9}
\sup_{t \in [0, \delta]} |E(Iu(t)) - E(Iu(0))| \leq O(\frac{1}{N^{3/2-}}).
\end{equation}
\end{lemma}

\noindent This implies global well-posedness for $u_{0} \in H^{s}(\mathbf{R}^{2})$, $s > 4/7$. Subsequent papers (see \cite{CKSTT3}, \cite{CR1}, \cite{CR}) have decreased the necessary regularity to

\begin{theorem}\label{t1.3}
$(\ref{1.1})$ has a global solution for $u_{0} \in H^{s}(\mathbf{R}^{2})$, $s > 1/3$.
\end{theorem}

\noindent This was proved by combining the I-method, a modified energy functional, and almost Morawetz estimates. The method will be described in more detail in the subsequent sections. In addition, the almost Morawetz estimates will be improved, thus improving Theorem $\ref{t1.3}$ to

\begin{theorem}\label{t1.4}
$(\ref{1.1})$ has a global solution for $u_{0} \in H^{s}(\mathbf{R}^{2})$, $s > 1/4$.
\end{theorem}

\noindent In $\S 2$ the modified energy functional of \cite{CKSTT3} will be recalled, as well as a modified local well-posedness theorem. In $\S 3$, the Morawetz inequality for $u(t,x)$ will be proved (originally proved in \cite{CR1}),

\begin{equation}\label{1.9}
\| u \|_{L_{t,x}^{4}([0, T] \times \mathbf{R}^{2})}^{4} \lesssim T^{1/3} \| u_{0} \|_{L^{2}(\mathbf{R}^{2})}^{3} \| u \|_{L_{t}^{\infty}([0, T], \dot{H}^{1}(\mathbf{R}^{2}))} + T^{1/3} \| u_{0} \|_{L^{2}(\mathbf{R}^{2})}^{4}.
\end{equation}

\noindent In $\S 4$, the known almost-Morawetz estimate in \cite{CR} for $Iu(t,x)$ will be improved. Finally, in $\S 5$, this improvement will be used to prove Theorem $\ref{t1.4}$.

\section{Modified Energy Functional}
In this section the known results concerning the modified energy functional will merely be stated. All of these results have been proved before (see \cite{CKSTT3} and \cite{CR}). If $u(t,x)$ solves $(\ref{1.1})$, then $Iu(t,x)$ solves

\begin{equation}\label{2.1}
i Iu_{t} + \Delta Iu = I(|u|^{2} u).
\end{equation}

\noindent If the nonlinearity was of the form $|Iu|^{2}(Iu)$, then $E(Iu(t))$ would be conserved. However, since $|Iu|^{2} Iu \neq I(|u|^{2} u)$,

\begin{equation}\label{2.2}
\partial_{t} E(Iu(t)) = 2 Re \int (Iu_{t}(t,x)) \{ I(|u(t,x)|^{2} \overline{u(t,x)}) - |Iu(t,x)|^{2} \overline{Iu(t,x)} \} dx.
\end{equation}

\noindent The change in energy decreases as $N \rightarrow \infty$.

\begin{theorem}\label{t2.1}
If $E(Iu(0)) \leq 1$, then there exists $\delta > 0$ such that

\begin{equation}\label{2.3}
|E(Iu(t)) - E(Iu(0))| \leq O(\frac{1}{N^{3/2-}}),
\end{equation}

\noindent for $t \in [0, \delta]$.
\end{theorem}

\noindent \emph{Proof:} See \cite{CKSTT1}.\vspace{5mm}

\noindent In \cite{CKSTT3}, the authors proved the existence of a modified energy functional $\tilde{E}(u(t))$ satisfying the properties:

1. $\tilde{E}(u(t))$ has a slower variation than $E(Iu(t))$.\vspace{5mm}

2. $\tilde{E}(u(t))$ is close to $E(Iu(t))$ in the sense that $E(Iu(t))$ can be controlled by $\tilde{E}(u(t))$.

\begin{proposition}\label{p2.1}
There exists a modified energy functional $\tilde{E}$ satisfying the fixed time estimate,

\begin{equation}\label{2.4}
|\tilde{E}(u(t)) - E(Iu(t))| \lesssim \frac{1}{\theta N^{2-}} \| Iu(t) \|_{H^{1}(\mathbf{R}^{2})}^{4}.
\end{equation}
\end{proposition}

\noindent \emph{Proof:} See $\S 4$ of \cite{CKSTT3}.

\begin{proposition}\label{p2.2}
\noindent $\tilde{E}(u(t))$ has the energy increment for a time interval J,

\begin{equation}\label{2.5}
|\sup_{t \in J} \tilde{E}(u(t)) - \tilde{E}(u(a))| \lesssim (\frac{\theta^{1/2}}{N^{3/2-}} + \frac{1}{N^{2-}} + \frac{1}{\theta N^{3-}}) \| Iu \|_{X^{1,1/2+}(J \times \mathbf{R}^{2})}^{4}.
\end{equation}
\end{proposition}

\noindent \noindent \emph{Proof:} See $\S 7$ and $\S 8$ of \cite{CKSTT3}.\vspace{5mm}

\noindent The $X^{1, 1/2+}$ norm will not be defined in this paper, because it will not be needed.

\begin{proposition}\label{p2.3}
Assume that

\begin{equation}\label{2.6}
\sup_{t \in J} E(Iu(t)) \leq 2,
\end{equation}

\noindent and for some $\epsilon > 0$,

\begin{equation}\label{2.7}
\| Iu \|_{L_{t,x}^{4}(J \times \mathbf{R}^{2})} \leq \epsilon.
\end{equation}

\noindent then

\begin{equation}\label{2.8}
|\sup_{t \in J} \tilde{E}(u(t)) - \tilde{E}(u(a))| \lesssim \frac{1}{N^{2-}} + \frac{\theta^{1/2}}{N^{3/2-}} + \frac{1}{\theta N^{3-}}.
\end{equation}

\noindent In particular, taking $\theta = \frac{1}{N}$ implies

\begin{equation}\label{2.9}
\sup_{t, t' \in J} |\tilde{E}(u(t)) - \tilde{E(u(t'))}| \lesssim \frac{1}{N^{2-}}.
\end{equation}
\end{proposition}

\noindent \emph{Proof:} See $\S 4$ of \cite{CR}.

\begin{theorem}\label{t2.4}
Let $$\| \langle \nabla \rangle Iu_{0} \|_{L^{2}(\mathbf{R}^{2})} = 1$$ and

\begin{equation}\label{2.10}
\int_{J_{k}} \int |Iu(t,x)|^{4} dx dt < \mu_{0},
\end{equation}

\noindent for some $\mu_{0} > 0$ sufficiently small. Then $(\ref{1.1})$ is locally well-posed on $[0, T]$ and

\begin{equation}\label{2.11}
Z_{I}(J_{k}, u) = \sup_{(q,r) \text{ admissible}} \| \langle \nabla \rangle Iu(t,x) \|_{L_{t}^{q} L_{x}^{r}(J_{k} \times \mathbf{R}^{2})} \leq C.
\end{equation}

\noindent $(q,r)$ is an admissible pair if $$\frac{2}{q} = 2(\frac{1}{2} - \frac{1}{r})$$ and $q > 2$.
\end{theorem}

\noindent \emph{Proof:} See $\S 3$ of \cite{CR1}.

\section{Morawetz inequalities}
In this section we will recall the proof of the following Morawetz inequality from \cite{CR1}. This recollection will be useful for the arguments given in the next section.

\begin{proposition}\label{p7.1}
If $u(t,x)$ solves $(\ref{1.1})$ then

\begin{equation}\label{7.0}
\| u(t,x) \|_{L_{t,x}^{4}([0, T] \times \mathbf{R}^{2})}^{4} \lesssim T^{1/3} \| u_{0} \|_{L^{2}(\mathbf{R}^{2})}^{3} \| u(t,x) \|_{L_{t}^{\infty}([0, T], \dot{H}^{1}(\mathbf{R}^{2}))} + T^{1/3} \| u_{0} \|_{L^{2}(\mathbf{R}^{2})}^{4}.
 \end{equation}
\end{proposition}

\noindent \emph{Proof:} Suppose that $v(t,z)$ solves the partial differential equation

\begin{equation}\label{7.1}
i v_{t} + \Delta_{z} v = F.
\end{equation}

\noindent Then define the quantities

\begin{equation}\label{7.2}
T_{0j}(t,z) = 2 Im(\overline{u(t,z)} \partial_{j} u(t,z)),
\end{equation}

\begin{equation}\label{7.3}
L_{jk}(t,z) = -\partial_{j} \partial_{k} (|u|^{2}) + 4 Re (\overline{\partial_{j} u(t,z)} \partial_{k} u(t,z)).
\end{equation}

\noindent These quantities obey the relation,

\begin{equation}\label{7.4}
\aligned
\partial_{t} T_{0j} + \partial_{k} L_{jk} =  2 (\overline{F(t,z)} \partial_{j} u(t,z) - \overline{u(t,z)} \partial_{j} F(t,z) \\ + F(t,z) \overline{\partial_{j} u(t,z)} - u(t,z) \overline{\partial_{j} F(t,z)}).
\endaligned
\end{equation}

\noindent Let $v(t,z)$ be a tensor product of solutions to $(\ref{1.1})$ on $\mathbf{R}^{2} \times \mathbf{R}^{2}$,

\begin{equation}\label{7.5}
(u_{1} \otimes u_{2})(t,z) = u_{1}(t,x) u_{2}(t,y) = v(t,z),
\end{equation}

\begin{equation}\label{7.6}
\aligned
i v_{t} + \Delta v = i \partial_{t} (u_{1}(t,x)) u_{2}(t,y) + i u_{1}(t,x) \partial_{t} (u_{2}(t,y)) \\ + (\Delta_{x} u_{1}(t,x)) u_{2}(t,y) + u_{1}(t,x) (\Delta_{y} u_{2}(t,y)) \\ = |u_{1}(t,x)|^{2} u_{1}(t,x) u_{2}(t,y) + |u_{2}(t,y)|^{2} u_{1}(t,x) u_{2}(t,y).
\endaligned
\end{equation}

\noindent Define the Morawetz action,

\begin{equation}\label{7.7}
M_{a}^{\otimes 2}(t) = 2 \int_{\mathbf{R}^{2} \times \mathbf{R}^{2}} \nabla a(z) \cdot Im(\overline{v(t,z)} \nabla v(t,z)) dz,
\end{equation}

\begin{equation}\label{7.8}
\partial_{t} M_{a}^{\otimes_{2}}(t) = 2 \int \partial_{j} a(z) \partial_{t} T_{0j}(t,z) dz,
\end{equation}

\noindent following the convention that repeated indices are summed.

\begin{equation}\label{7.9.1}
\partial_{t} M_{a}^{\otimes_{2}}(t) = 2 \int \partial_{j} \partial_{kk} (|v|^{2}) \partial_{j} a(z)
\end{equation}

\begin{equation}\label{7.9.2}
- 8 \int \partial_{k} Re(\overline{\partial_{j} v(t,z)} \partial_{k} v(t,z)) \partial_{j} a(z) dz
\end{equation}

\begin{equation}\label{7.9.3}
\aligned
+ 4 \int \{ \overline{F(t,z)} \partial_{j} v(t,z) - \overline{v(t,z)} \partial_{j} F(t,z)\\ + F(t,z) \overline{\partial_{j} v(t,z)} - v(t,z) \overline{\partial_{j} F(t,z)} \} \partial_{j} a(z).
\endaligned
\end{equation}

\vspace{5mm}

\noindent Let $v(t,z) = u(t,x) u(t,y)$, where $u$ solves $(\ref{1.1})$. Take the term $(\ref{7.9.1})$ first.

$$2 \int \partial_{j} \partial_{kk} (|v(t,z)|^{2}) \partial_{j} a(z) = -2 \int |v(t,z)|^{2} (\Delta \Delta a(z)) dz.$$

\noindent Now let $a(z) = a(x,y) = f(|x - y|)$, where $f$ is a smooth, convex function. Let

\begin{equation}\label{7.11}
f(x) = \left\{
         \begin{array}{ll}
           \frac{1}{2M} x^{2}(1 - \log \frac{x}{M}), & \hbox{if $|x| < \frac{M}{\sqrt{e}}$;} \\
           100x, & \hbox{if $|x| > M$.}
         \end{array}
       \right.
\end{equation}

\noindent For $|x - y| < \frac{M}{\sqrt{e}}$, $$\Delta a(x,y) = \frac{2}{M} \log(\frac{M}{|x - y|}) \Rightarrow -\Delta \Delta a(x,y) = \frac{2}{M} \delta_{x = y},$$ and for $|x - y| > M$, $$-\Delta \Delta a(x,y) = O(\frac{1}{|x - y|^{3}}) = O(\frac{1}{M^{3}}).$$

$$\int_{0}^{T} \int_{\mathbf{R}^{2} \times \mathbf{R}^{2}} (-\Delta \Delta a(x,y)) |u(t,x)|^{2} |u(t,y)|^{2} dx dy dt = \frac{2}{M} \int_{0}^{T} \int_{\mathbf{R}^{2}} |u(t,x)|^{4} dx dt$$

$$ + O(\frac{1}{M^{3}}) \int_{0}^{T} \int_{\mathbf{R}^{2} \times \mathbf{R}^{2}} |u(t,x)|^{2} |u(t,y)|^{2} dx dy dt.$$

\noindent Since $M$ will be large, $|\nabla a(z)|$ is uniformly bounded on $\mathbf{R}^{2} \times \mathbf{R}^{2}$, and

$$|M_{a}^{\otimes 2}(t)| = 2 |\int_{\mathbf{R}^{2} \times \mathbf{R}^{2}} \nabla a(z) \cdot Im(\overline{v(t,z)} \nabla v(t,z)) dz|$$

$$ \lesssim  \| u_{1}(t,x) \|_{L_{t}^{\infty}([0, T], L^{2}(\mathbf{R}^{2}))}^{2} \| u_{2}(t,y) \|_{L_{t}^{\infty}([0, T], \dot{H}^{1}(\mathbf{R}^{2}))} \| u_{2}(t,y) \|_{L_{t}^{\infty}([0, T], L^{2}(\mathbf{R}^{2}))}$$

$$ + \| u_{2}(t,y) \|_{L_{t}^{\infty}([0, T], L^{2}(\mathbf{R}^{2}))}^{2} \| u_{1}(t,x) \|_{L_{t}^{\infty}([0, T], \dot{H}^{1}(\mathbf{R}^{2}))} \| u_{1}(t,x) \|_{L_{t}^{\infty}([0, T], L^{2}(\mathbf{R}^{2}))}.$$

\noindent This implies,

\begin{equation}\label{7.12}
\aligned
\frac{2}{M} \int_{0}^{T} \int_{\mathbf{R}^{2}} |u(t,x)|^{4} dx dt + O(\frac{1}{M^{3}}) \int_{0}^{T} \int_{\mathbf{R}^{2} \times \mathbf{R}^{2}} |u(t,x)|^{2} |u(t,y)|^{2} dx dy dt \\
+ (\ref{7.9.2}) + (\ref{7.9.3}) \lesssim \| u_{0} \|_{L^{2}(\mathbf{R}^{2})}^{3} \| u \|_{L_{t}^{\infty} \dot{H}_{x}^{1}([0, T] \times \mathbf{R}^{3})}.
\endaligned
\end{equation}

\noindent The proof will be complete once we prove $(\ref{7.9.2})$ and $(\ref{7.9.3})$ are positive.

\begin{lemma}\label{l7.2}
Let $f$ be a convex function. Then $$\partial_{j} \partial_{k} a(z),$$ gives a positive definite matrix for all $z \in \mathbf{R}^{2} \times \mathbf{R}^{2}$ if $a(z) = f(|x - y|)$.
\end{lemma}

\noindent \emph{Proof:} $$\partial_{j} \partial_{k} f(|x - y|) = f''(|x - y|) \frac{(x - y)_{j} (x - y)_{k}}{|x - y|^{2}} + \frac{f'(|x - y|)}{|x - y|}(\delta_{jk} - \frac{(x - y)_{j} (x - y)_{k}}{|x - y|^{2}}).$$

\noindent Take the inner product defined by this matrix.

$$\langle z_{j} z_{k} | f''(| x - y|) \frac{(x - y)_{j} (x - y)_{k}}{|x - y|^{2}} \rangle = \frac{f''(|x - y|)}{|x - y|^{2}} (z \cdot (x - y))^{2}.$$

$$|\langle z_{j} z_{k} | \frac{f'(|x - y|)}{|x - y|} \frac{(x - y)_{j} (x - y)_{k}}{|x - y|^{2}} \rangle| \leq \frac{|f'(|x - y|)|}{|x - y|} |z|^{2},$$

$$\langle z_{j} z_{k} | \frac{f'(|x - y|)}{|x - y|} \delta_{jk} \rangle = \frac{f'(|x - y|)}{|x - y|} |z|^{2}.$$

\noindent This proves the lemma. $\Box$\vspace{5mm}

\noindent In particular, after integrating by parts, $(\ref{7.9.2}) \geq 0$.\vspace{5mm}

\noindent To evaluate $(\ref{7.9.3})$, without loss of generality take $j = 1$.

$$\overline{F(t,z)} \partial_{1} v(t,z) - \overline{v(t,z)} \partial_{1} F(t,z)$$

$$ = |u(t,y)|^{2} \overline{u(t,x) u(t,y)} \partial_{1}(u(t,x) u(t,y)) - \overline{u(t,x) u(t,y)} \partial_{1}(|u(t,y)|^{2} u(t,y) u(t,x))$$

$$ + |u(t,x)|^{2} \overline{u(t,x) u(t,y)} \partial_{1}(u(t,x) u(t,y)) - \overline{u(t,x) u(t,y)} \partial_{1}(|u(t,x)|^{2} u(t,y) u(t,x)).$$

$$|u(t,y)|^{2} \overline{u(t,x) u(t,y)} \partial_{1}(u(t,x) u(t,y)) - \overline{u(t,x) u(t,y)} \partial_{1}(|u(t,y)|^{2} u(t,y) u(t,x)) = 0.$$

$$|u(t,x)|^{2} \overline{u(t,x) u(t,y)} \partial_{1}(u(t,x) u(t,y)) - \overline{u(t,x) u(t,y)} \partial_{1}(|u(t,x)|^{2} u(t,y) u(t,x))$$

$$ = -|u(t,x)|^{2} \overline{u(t,x) u(t,y)} \partial_{1}(u(t,x) u(t,y)) - |u(t,x)|^{2} (u(t,x) u(t,y)) \partial_{1}(\overline{u(t,x) u(t,y)}$$

$$ = \frac{-1}{2} \partial_{1}(|u(t,x)|^{4} |u(t,y)|^{2}).$$

\noindent Similarly,

\begin{equation}\label{7.13}
F \overline{\partial_{1} u} - u \overline{\partial_{1} F} = \frac{-1}{2} \partial_{1}(|u(t,x)|^{4} |u(t,y)|^{2}).
\end{equation}

\noindent Make a similar calculation $j = 2, 3, 4$, although when $j = 3$ or $4$ switch $x$ and $y$ in $(\ref{7.13})$. Therefore, $(\ref{7.9.3})$ is a sum of terms of the form

$$-\int_{0}^{T} \int_{\mathbf{R}^{2}} \int_{\mathbf{R}^{2}} \partial_{j}(|u(t,x)|^{4} |u(t,y)|^{2}) a_{j}(z) dx dy dt,$$

\noindent when $j = 1, 2$ and

$$-\int_{0}^{T} \int_{\mathbf{R}^{2}} \int_{\mathbf{R}^{2}} \partial_{j}(|u(t,y)|^{4} |u(t,x)|^{2}) a_{j}(z) dx dy dt,$$

\noindent when $j = 3, 4$. Integrating by parts and noticing $$a_{jj}(z) = f''(|x - y|) \frac{(x - y)_{j}^{2}}{|x - y|^{2}} + \frac{f'(|x - y|)}{|x - y|} (1 - \frac{(x - y)_{j}^{2}}{|x - y|^{2}}) \geq 0$$

\noindent proves $(\ref{7.9.3}) \geq 0$. Combining terms,

$$\frac{2}{M} \int_{0}^{T} \int_{\mathbf{R}^{2}} |u(t,x)|^{4} dx dt \lesssim \sup_{[0, T]} \| u(t,x) \|_{2}^{3} \| u(t,x) \|_{\dot{H}^{1}} + O(\frac{T}{M^{3}}) \sup_{[0, T]} \| u(t,x) \|_{2}^{4}.$$

\noindent Choosing $M = T^{1/3}$ proves the proposition. $\Box$

\section{Almost Morawetz Inequalities}
\noindent In this section, the almost Morawetz estimate in \cite{CR1}, \cite{CR} will be improved. For $u_{0}$ with regularity below $s = 1$, if $u(t,x)$ solves $(\ref{1.1})$ then $Iu(t,x)$ solves

\begin{equation}\label{8.1}
i Iu(t,x) + \Delta Iu(t,x) = I(|u(t,x)|^{2} u(t,x).
\end{equation}

\begin{proposition}\label{p8.1}
\noindent Define the quantity

\begin{equation}\label{8.2}
Z_{I}([0, T]) = \sup_{(q,r) \text{ admissible}} \| \langle D \rangle Iu \|_{L_{t}^{q} L_{x}^{r}([0, T] \times \mathbf{R}^{2})}.
\end{equation}

\begin{equation}\label{8.3}
\aligned
\| Iu(t,x) \|_{L_{t,x}^{4}([0, T] \times \mathbf{R}^{2})}^{4} \lesssim T^{1/3} \| u_{0} \|_{L^{2}(\mathbf{R}^{2})}^{2} \| Iu(t,x) \|_{L_{t}^{\infty}([0, T]; \dot{H}^{1}(\mathbf{R}^{2}))} \\ + T^{1/3} \| u_{0} \|_{L^{2}(\mathbf{R}^{2})}^{4} + T^{1/3} \sum_{k} \frac{Z_{I}(J_{k})^{6}}{N^{2-}},
\endaligned
\end{equation}

\noindent where $J_{k}$ is a partition of $[0, T]$.
\end{proposition}

\noindent \emph{Proof:} Split the nonlinearity

\begin{equation}\label{8.4}
\aligned
F &= I(|u(t,x)|^{2} u(t,x)) Iu(t,y) + Iu(t,x) I(|u(t,y)|^{2} u(t,y)) = \mathcal N_{g} + \mathcal N_{b}, \\
\mathcal N_{g} &= |Iu(t,x)|^{2} Iu(t,x) Iu(t,y) + |Iu(t,y)|^{2} Iu(t,x) Iu(t,y), \\
\mathcal N_{b} &= F - \mathcal N_{g}.
\endaligned
\end{equation}

\noindent After taking a tensor product of solutions $v(t,z) = Iu(t,x) Iu(t,y)$, repeat the procedure from $\S 3$ to obtain

\begin{equation}\label{8.5}
\aligned
-2 \int_{0}^{T} |v(t,z)|^{2} (\Delta \Delta a(z)) dz \\
+ 8 \int Re(\overline{\partial_{j} v(t,z)} \partial_{k} v(t,z)) dz \\
+ 4 \int (\overline{F(t,z)} \partial_{j} v(t,z) - \overline{v(t,z)} \partial_{j} F(t,z) \\ + F(t,z) \overline{\partial_{j} v(t,z)} - v(t,z) \overline{\partial_{j} F(t,z)}) \partial_{j} a(z) dz \\= M_{a}(T) - M_{a}(0)
\endaligned
\end{equation}

\noindent Once again, the second term $8 \int Re(\overline{\partial_{j} v(t,z)} \partial_{k} v(t,z)) dz$ is strictly positive and can be discarded, as well as the parts of the third term with $\mathcal N_{g}$ in place of F. Therefore

\begin{equation}\label{8.6}
\aligned
\int_{0}^{T} \int |Iu(t,x)|^{4} dx dt \lesssim T^{1/3} \| u_{0} \|_{2}^{3} \| Iu \|_{L_{t}^{\infty}([0, T], \dot{H}^{1}(\mathbf{R}^{2}))} + T^{1/3} \| u_{0} \|_{2}^{4} \\
+ T^{1/3} \int_{0}^{T} \int (\overline{\mathcal N_{b}} \partial_{j} v(t,z) - \overline{v(t,z)} \partial_{j} \mathcal N_{b} + \mathcal N_{b} \overline{\partial_{j} v(t,z)} - v(t,z) \overline{\partial_{j} \mathcal N_{b}}) \partial_{j} a(z) dz
\endaligned
\end{equation}

\noindent To handle $$T^{1/3} \int_{0}^{T} \int (\overline{\mathcal N_{b}} \partial_{j} v(t,z) - \overline{v(t,z)} \partial_{j} \mathcal N_{b} + \mathcal N_{b} \overline{\partial_{j} v(t,z)} - v(t,z) \overline{\partial_{j} \mathcal N_{b}}) \partial_{j} a(z) dz$$ it suffices to handle terms of the form

\begin{equation}\label{8.7}
\int_{J_{k}} \int_{\mathbf{R}^{2} \times \mathbf{R}^{2}} \nabla a \cdot (\mathcal N_{b}) \overline{\nabla v(t,x)} dz dt,
\end{equation}

\noindent as well as terms of the form

\begin{equation}\label{8.8}
\int_{J_{k}} \int_{\mathbf{R}^{2} \times \mathbf{R}^{2}} \nabla a \cdot (\nabla \mathcal N_{b}) \overline{v(t,z)} dz dt.
\end{equation}

\noindent Integrating by parts in $x$, $(\ref{8.8})$ is a sum of terms of the form $(\ref{8.7})$, along with terms of the form

\begin{equation}\label{8.9}
\int_{0}^{T} \int \Delta a(z) \overline{\mathcal N_{b}(t,z)} v(t,z) dz dt.
\end{equation}

\noindent $(\ref{8.7})$ will be tackled first.

$$\int_{J_{k}} \int_{\mathbf{R}^{2} \times \mathbf{R}^{2}} \nabla a \cdot (\mathcal N_{b}) \overline{\nabla v(t,z)} dz dt$$

\noindent is a sum of terms of the form

\begin{equation}\label{8.10}
\int_{J_{k}} \int \mathcal N_{b}(t,z) |Iu(t,y)| |\nabla Iu(t,x)| dx dy dt.
\end{equation}

$$\mathcal N_{b} = Iu(t,x)[I(|u(t,y)|^{2} u(t,y)) - |Iu(t,y)|^{2} Iu(t,y)]$$ $$+ Iu(t,y)[I(|u(t,x)|^{2} u(t,x)) - |Iu(t,x)|^{2} Iu(t,x)].$$

\noindent This implies

\begin{equation}\label{8.11}
(\ref{8.10}) \lesssim \| I(|u(t,x)|^{2} u(t,x)) - |Iu(t,x)|^{2} Iu(t,x) \|_{L_{t}^{1} L_{x}^{2}(J_{k} \times \mathbf{R}^{2})} Z_{I}(J_{k})^{3}.
\end{equation}

\noindent The quantity

\begin{equation}\label{8.12}
\| I(|u(t,x)|^{2} u(t,x)) - |Iu(t,x)|^{2} Iu(t,x) \|_{L_{t}^{1} L_{x}^{2}}
\end{equation}

\noindent can be estimated by making a Littlewood-Paley partition of $u(t,x)$. Define a quantity $F(t,\xi)$

$$F(t,\xi) = \int_{\xi_{1} + \xi_{2} + \xi_{3} = \xi} [m(\xi_{1} + \xi_{2} + \xi_{3}) - m(\xi_{1}) m(\xi_{2}) m(\xi_{3})]$$ $$\times  \hat{u}(t,\xi_{1}) \hat{\bar{u}}(t,\xi_{2}) \hat{u}(t,\xi_{3}) d\xi_{1} d\xi_{2}$$

$$ = \int_{\xi_{1}+ \xi_{2} + \xi_{3} = \xi} \frac{[m(\xi) - m(\xi_{1}) m(\xi_{2}) m(\xi_{3})]}{m(\xi_{1}) m(\xi_{2}) m(\xi_{3})} \widehat{Iu}(t,\xi_{1}) \widehat{\overline{Iu}}(t,\xi_{2}) \hat{Iu}(t,\xi_{3}).$$

\noindent Suppose $\hat{u}(t,\xi_{i})$ is supported on the frequency region $|\xi_{i}| \sim N_{i}$, and without loss of generality suppose $N_{1} \geq N_{2} \geq N_{3}$. Consider four regions separately.\vspace{5mm}

\noindent \emph{$N_{1} << N$:} In this case the multipliers $m(\xi_{i}) = 1$, so $$\frac{[m(\xi) - m(\xi_{1}) m(\xi_{2}) m(\xi_{3})]}{m(\xi_{1}) m(\xi_{2}) m(\xi_{3})} = 0.$$\vspace{5mm}

\noindent \emph{$N_{2} << N \lesssim N_{1}$:} By the fundamental theorem of calculus, $$m(\xi_{1} + \xi_{2} + \xi_{3}) - m(\xi_{1}) \lesssim |\xi_{2} + \xi_{3}| \nabla m(\xi_{1}).$$ $$\frac{\nabla m(\xi_{1})}{m(\xi_{1})} \lesssim \frac{1}{|\xi_{1}|}.$$

$$\| \frac{|\xi_{2} + \xi_{3}| |\xi_{1}|}{|\xi_{1}|^{2}} \widehat{Iu}_{1}(t,\xi_{1}) \widehat{\overline{Iu}}_{2}(t,\xi_{2}) \widehat{Iu}_{3}(t,\xi_{3}) \|_{L_{t}^{1} L_{x}^{2}}$$

$$\lesssim \frac{1}{N^{2}} \| \nabla Iu_{1} \|_{L_{t}^{3} L_{x}^{6}} \| \nabla Iu_{2} \|_{L_{t}^{3} L_{x}^{6}} \| Iu_{3} \|_{L_{t}^{3} L_{x}^{6}} \lesssim \frac{N_{1}^{-} N_{2}^{-} N_{3}^{-}}{N^{2-}} Z_{I}^{3}.$$

\noindent \emph{$N_{3} << N \lesssim N_{2} \leq N_{1}$:} In this case, make the trivial multiplier estimate,

$$|\frac{m(\xi_{1} + \xi_{2} + \xi_{3}) - m(\xi_{1}) m(\xi_{2})}{m(\xi_{1}) m(\xi_{2})}| \leq \frac{m(\xi_{1} + \xi_{2} + \xi_{3})}{m(\xi_{1}) m(\xi_{2})} + 1$$

$$\lesssim \frac{m(\xi_{1} + \xi_{2} + \xi_{3})}{m(\xi_{1}) m(\xi_{2})},$$

\noindent since $m(\xi_{1} + \xi_{2} + \xi_{3}) \sim m(\xi_{1})$ and $m(\xi_{2}) \leq 1$,

$$1 + \frac{m(\xi_{1} + \xi_{2} + \xi_{3})}{m(\xi_{1}) m(\xi_{2})} \lesssim \frac{1}{m(\xi_{2})}.$$

$$\frac{1}{m(\xi_{2}) |\xi_{2}| |\xi_{1}|} \lesssim \frac{1}{m(N) N |\xi_{1}|} \lesssim \frac{1}{N^{2}}.$$

\noindent This uses the fact that $m(\xi) \xi$ is monotone increasing for any $s > 0$ and $m(N) N = N$. Therefore,

$$\| \int_{\xi_{1} + \xi_{2} + \xi_{3} = \xi} \frac{m(\xi_{1} + \xi_{2} + \xi_{3})}{m(\xi_{1}) m(\xi_{2})} \widehat{Iu}(t,\xi_{1}) \widehat{\overline{Iu}}(t,\xi_{2}) \widehat{Iu}(t,\xi_{3}) d\xi_{1} d\xi_{2} \|_{L_{t}^{1} L_{x}^{2}} \lesssim$$

$$\| \int_{\xi_{1} + \xi_{2} + \xi_{3} = \xi} \frac{m(\xi_{1} + \xi_{2} + \xi_{3})}{m(\xi_{1}) m(\xi_{2}) |\xi_{1}| |\xi_{2}|} \widehat{\nabla Iu}(t,\xi_{1}) \widehat{\overline{\nabla Iu}}(t,\xi_{2}) \widehat{Iu}(t,\xi_{3}) d\xi_{1} d\xi_{2} \|_{L_{t}^{1} L_{x}^{2}}$$  $$\lesssim \frac{N_{1}^{-} N_{2}^{-} N_{3}^{-}}{N^{2-}} Z_{I}^{3}.$$

\noindent Finally, consider the region\vspace{5mm}

\noindent \emph{$N \lesssim N_{3} \leq N_{2} \leq N_{1}$:} Doing the same analysis,

$$\frac{m(\xi_{1} + \xi_{2} + \xi_{3})}{|\xi_{1}| |\xi_{2}| |\xi_{3}|} \frac{1}{m(\xi_{1}) m(\xi_{2}) m(\xi_{3})} \lesssim \frac{1}{N^{3}} Z_{I}^{3}.$$

\noindent This is because $m(\xi_{1} + \xi_{2} + \xi_{3}) \sim m(\xi_{1})$, and $m(\xi_{2}) |\xi_{2}| \gtrsim m(N) N$, $m(\xi_{3}) |\xi_{3}| \gtrsim m(N) N$.

$$\| \int_{\xi_{1} + \xi_{2} + \xi_{3} = \xi} \frac{m(\xi_{1} + \xi_{2} + \xi_{3})}{m(\xi_{1}) m(\xi_{2}) |\xi_{1}| |\xi_{2}| |\xi_{3}|} \widehat{\nabla Iu}(t,\xi_{1}) \widehat{\overline{\nabla Iu}}(t,\xi_{2}) \widehat{\nabla Iu}(t,\xi_{3}) d\xi \|_{L_{t}^{1} L_{x}^{2}}$$  $$\lesssim \frac{N_{1}^{-} N_{2}^{-} N_{3}^{-}}{N^{3-}} Z_{I}^{3}.$$

\noindent This proves the proposition for terms of the form $(\ref{8.7})$.\vspace{5mm}

\noindent Turning to $(\ref{8.9})$,

$$(\ref{8.9}) \lesssim \int_{J_{k}} \int | I(|u(t,x)|^{2} u(t,x)) - |Iu(t,x)|^{2}(Iu(t,x)) | Iu(t,x) \Delta a(|x - y|) |Iu(t,y)|^{2} dx dy dt.$$ On $|x - y| < \frac{M}{\sqrt{e}}$, $$\Delta a(x,y) = \frac{2}{M} \log(\frac{M}{|x - y|}),$$ and for large $|x - y|$, $$\Delta a(x,y) = O(\frac{1}{|x - y|}).$$ Therefore, for $|x - y| > 1$, $\Delta a(x,y)$ is uniformly bounded. This bound is uniform for $M \geq 1$.

$$
\int_{0}^{T} \int \int_{|x - y| > 1} |I(|u(t,x)|^{2} u(t,x)) - |Iu(t,x)|^{2} Iu(t,x)|$$ $$\times |Iu(t,x)| |Iu(t,y)|^{2} \Delta a(x,y) dx dy dt$$

$$ \leq \sup_{x} (\int_{|x - y| > 1} \Delta a(x,y) |Iu(t,y)|^{2} dy)$$ $$\times \int_{0}^{T} \int |I(|u(t,x)|^{2} u(t,x)) - |Iu(t,x)|^{2} Iu(t,x)| |Iu(t,x)| dx dt.$$

$$\int_{|x - y| > 1} \Delta a(x,y) |Iu(t,y)|^{2} \lesssim \| Iu(t,y) \|_{L^{2}}^{2}.$$

\begin{equation}\label{8.14}
\aligned
\int_{0}^{T} \int \int_{|x - y| > 1} |I(|u(t,x)|^{2} u(t,x)) - |Iu(t,x)|^{2} Iu(t,x)| \\ \times |Iu(t,x)| |Iu(t,y)|^{2} \Delta a(x,y) dx dy dt \\ \lesssim \| u_{0} \|_{2}^{3} \| I(|u(t,x)|^{2} u(t,x)) - |Iu(t,x)|^{2} Iu(t,x) \|_{L_{t}^{1} L_{x}^{2}(J_{k} \times \mathbf{R}^{2})}.
\endaligned
\end{equation}

\noindent For a fixed $x$ take the region $|x - y| \leq 1$,

$$\int_{|x - y| \leq 1} \Delta a(x,y) |Iu(t,y)|^{2} dy \leq \| Iu(t,y) \|_{L^{4}}^{2} \| \frac{2}{M} \log(\frac{M}{|x - y|}) \|_{L^{4}(|x - y| \leq 1)}.$$ Since $\| Iu(t,y) \|_{L^{4}(\mathbf{R}^{2})} \leq \| \langle \nabla \rangle^{1/2} Iu \|_{L^{2}(\mathbf{R}^{2})}$, therefore,

\begin{equation}\label{8.15}
\sup_{x} \int \Delta a(x,y) |Iu(t,y)|^{2} dy \leq C.
\end{equation}

$$\int_{J_{k}} \int \int |I(|u(t,x)|^{2} u(t,x)) - |Iu(t,x)|^{2} Iu(t,x)| |Iu(t,x)| |\Delta a(x,y)| |Iu(t,y)|^{2} dx dy dt$$ $$\lesssim \| I(|u(t,x)|^{2} u(t,x)) - |Iu(t,x)|^{2} Iu(t,x) \|_{L_{t}^{1} L_{x}^{2}(J_{k} \times \mathbf{R}^{2})} Z_{I}(J_{k})^{3} \lesssim \frac{1}{N^{2-}} Z_{I}(J_{k})^{6}.$$

\noindent This completes the proof of the proposition. $\Box$\vspace{5mm}

\section{Proof of Theorem $\ref{t1.4}$}
Fix a time interval $[0, T_{0}]$. We wish to show that $(\ref{1.1})$ has a solution on that time interval. If $u(t,x)$ is a solution on $[0, T]$ then $$\frac{1}{\lambda} u(\frac{t}{\lambda^{2}}, \frac{x}{\lambda})$$ is a solution on $[0, \lambda^{2} T]$. Let $u_{0,\lambda}$ denote the rescaled solution at $t = 0$, and let $u_{\lambda}(t)$ be the rescaled solution. $$\| u_{0, \lambda} \|_{\dot{H}^{s}(\mathbf{R}^{2})} = \lambda^{-s} \| u_{0} \|_{\dot{H}^{s}(\mathbf{R}^{2})}.$$

\begin{equation}\label{5.1}
\| Iu_{0} \|_{\dot{H}^{1}} \lesssim N^{1 - s} \| u_{0} \|_{\dot{H}^{s}(\mathbf{R}^{2})}.
\end{equation}

\noindent Choose $\lambda = C_{0}(\| u_{0} \|_{H^{s}(\mathbf{R}^{2})}) N^{(1 - s)/s}$ so that $$E(Iu_{0, \lambda}) = \frac{2}{5}.$$ We now wish to prove $E(Iu_{\lambda}(t)) \leq 1$ on $[0, \lambda^{2} T_{0}]$.\vspace{5mm}

\noindent Next, define a subset of $[0, \lambda^{2} T_{0}]$,

\begin{equation}\label{5.2}
F_{T} = \{ t : \tilde{E}(u_{\lambda}(t)) \leq \frac{3}{4} \}.
\end{equation}

\noindent By the fixed time estimate $(\ref{2.4})$, $|\tilde{E}_{\lambda}(u(0)) - E_{\lambda}(Iu(0))| \lesssim \frac{1}{\theta N^{2-}}$,  assume $\tilde{E}(u_{\lambda}(0)) \leq \frac{1}{2}$ since $E(Iu_{\lambda}(0)) \leq \frac{2}{5}$, therefore, $0 \in F_{T}$. Furthermore, $F_{T}$ is closed in $[0, \lambda^{2} T_{0}]$ by the dominated convergence theorem. It remains therefore to show $F_{T}$ is open in $[0, \lambda^{2} T_{0}]$. If $\tilde{E}(u_{\lambda}(t)) \leq \frac{3}{4}$ on $[0, T']$, then for some $\delta > 0$, $\tilde{E}(u_{\lambda}(t)) \leq \frac{9}{10}$ on $[0, T' + \delta]$, which in turn implies $E(Iu_{\lambda}(t)) \leq 1$ on $[0, T' + \delta]$.\vspace{5mm}

\noindent Because $E(Iu_{\lambda}(t)) \leq 1$ on $[0, T' + \delta]$, by the Sobolev embedding theorem $\| Iu_{\lambda}(t,x) \|_{L_{t,x}^{4}([0, T' + \delta] \times \mathbf{R}^{2})}$ is finite. Next apply the local well-posedness theorem $\ref{t2.4}$. If $\| Iu_{\lambda} \|_{L_{t,x}^{4}(J_{k} \times \mathbf{R}^{2})} \leq \epsilon$ and $\| \langle \nabla \rangle Iu_{0, \lambda} \|_{L^{2}(\mathbf{R}^{2})} \leq 1$, then

\begin{equation}\label{5.3}
Z(J_{k}, u_{\lambda}) \leq C.
\end{equation}

\noindent The interval $[0, T' + \delta]$ can be partitioned into $$\frac{\| Iu_{\lambda} \|_{L_{t,x}^{4}([0, T' + \delta] \times \mathbf{R}^{2})}^{4}}{\epsilon^{4}}$$ pieces $J_{k}$ such that $Z(J_{k}) \leq C$. Next, apply the almost Morawetz inequality.

\begin{equation}\label{5.4}
\aligned
\| Iu_{\lambda} \|_{L_{t,x}^{4}([0, T' + \delta])}^{4} &\leq C' \lambda^{2/3} T_{0}^{1/3}(\| u_{0} \|_{L^{2}}^{4} + \| u_{0} \|_{L^{2}}^{3} \| Iu_{\lambda} \|_{L_{t}^{\infty} \dot{H}_{x}^{1}(\mathbf{R}^{2})} + \frac{1}{N^{2-}} \sum_{k} Z_{k}(J_{k}, u_{\lambda})^{6}) \\
&\leq C' \lambda^{2/3} T_{0}^{1/3}(\| u_{0} \|_{L^{2}}^{4} + \| u_{0} \|_{L^{2}}^{3} \sup_{[0, T' + \delta]} E(Iu_{\lambda}(t)) + \frac{C^{6}}{\epsilon^{4} N^{2-}} \| Iu_{\lambda} \|_{L_{t,x}^{4}}^{4}).
\endaligned
\end{equation}

$$\lambda^{2/3} N^{-2+} \sim N^{\frac{2 - 8s}{3s}+},$$ so for $s > 1/4$, choosing $N$ sufficiently large, $$\frac{C' C^{6}}{\epsilon^{4}} \lambda^{2/3} T_{0}^{1/3} \leq \frac{1}{2}.$$ Therefore, the remainder can be absorbed into the left hand side and

$$\| Iu_{\lambda} \|_{L_{t,x}^{4}([0, T' + \delta])}^{4} \leq 2 C' \lambda^{2/3} T_{0}^{1/3} (\| u_{0} \|_{L^{2}}^{3} + \| u_{0} \|_{L^{2}}^{4}).$$

\noindent Partitioning $[0, T' + \delta]$ into $\frac{2 C_{0}}{\epsilon^{4}} \lambda^{2/3} T^{1/3} (\| u_{0} \|_{L^{2}}^{3} + \| u_{0} \|_{L^{2}}^{4})$ pieces,

\begin{equation}\label{5.5}
\sup_{t \in [0, T' + \delta]} |\tilde{E}(u_{\lambda}(t))| \leq \frac{1}{2} + \frac{2 C' \lambda^{2/3} T_{0}^{1/3} (\| u_{0} \|_{2}^{4} + \| u_{0} \|_{2}^{3})}{\epsilon^{4} N^{2-}}.
\end{equation}

\noindent Taking $N(T_{0}, \| u_{0} \|_{2})$ sufficiently large, this implies $F_{T} = [0, \lambda^{2} T_{0}]$.\vspace{5mm}

\noindent \emph{Proof of Theorem $\ref{t1.4}$:} Let $$N = (\frac{20 C' C_{0}^{2/3} T_{0}^{1/3} (m_{0}^{4} + m_{0}^{3})}{\epsilon^{4}})^{3s/(8s - 2)+} = C(m_{0}) T_{0}^{\frac{s}{8s - 2}+}.$$

\noindent This implies $$\sup_{[0, \lambda^{2} T_{0}]} \tilde{E}(u_{\lambda}(t)) \leq \frac{3}{5},$$ which in turn implies $E(Iu_{\lambda}(t)) \leq 1$ on $[0, \lambda^{2} T_{0}]$. Splitting the solution $u_{\lambda}(t) = P_{\leq N} u_{\lambda}(t) + P_{> N} u_{\lambda}(t)$, $$\| P_{\leq N} u_{\lambda} \|_{H^{s}(\mathbf{R}^{2})} \leq \| I_{N} u_{\lambda} \|_{H^{1}(\mathbf{R}^{2})} \leq E(Iu_{\lambda}(t)),$$ $$\| P_{> N} u_{\lambda} \|_{H^{s}(\mathbf{R}^{2})} \leq N^{s - 1} \| I_{N} u_{\lambda} \|_{H^{1}(\mathbf{R}^{2})} \leq N^{s - 1} E(Iu_{\lambda}(t)),$$ which proves that $$\| u_{\lambda}(t) \|_{H^{s}(\mathbf{R}^{2})} \leq 2.$$ Finally, $\lambda = C_{0} N^{\frac{1 - s}{s}} = C(m_{0}) C_{0}(\| u_{0} \|_{H^{s}}) T_{0}^{\frac{1 - s}{8s - 2}+}$, so rescaling back,

\begin{equation}\label{5.6}
\sup_{t \in [0, T]} \| u(t) \|_{H^{s}(\mathbf{R}^{2})} \leq C(m_{0}) C_{0}(\| u_{0} \|_{H^{s}(\mathbf{R}^{2})}) T^{\frac{s(1 - s)}{8s - 2}+}.
\end{equation}

\noindent This proves the theorem. $\Box$

\newpage
\nocite*
\bibliographystyle{plain}
\bibliography{improvement}
\end{document}